\documentclass[11pt]{amsart}
\usepackage{amsmath,amsfonts,amscd,amssymb}
\theoremstyle{plain}
\newcounter{thmcount}
\newtheorem{theorem}[thmcount]{Theorem}
\newtheorem{proposition}[thmcount]{Proposition}
\newtheorem{lemma}[thmcount]{Lemma}

\theoremstyle{definition}
\newtheorem{remark}[thmcount]{Remark}

\def\O{{\mathcal O}}
\def\F{{\mathbb F}}

\def\Q{{\mathbb Q}}

\def\C{{\mathbb C}}
\def\Z{{\mathbb Z}}

\def\D{{\mathfrak D}}
\let\iff\Leftrightarrow
\def\triv#1{{\mathbf 1}_{#1}}
\def\V{V}

\def\newmathop#1{\expandafter\gdef\csname #1\endcsname{\mathop{\rm #1}\nolimits}}
\newmathop{tr}
\newmathop{rk}
\newmathop{ord}
\newmathop{Gal}
\newmathop{Ind}
\newmathop{Res}
\newmathop{Frob}
\newmathop{GL}
\newmathop{SL}
\newmathop{Aut}

\def\notsubset{\not\subset}

\def\Cy#1{C_{#1}}
\def\Di#1{D_{#1}}
\def\Qu{Q_8}
\def\Sy{H_{16}}
\def\Bo{B_{12}}
\def\Sl{\SL_{24}}
\def\Gl{\GL_{48}}
\def\Vi{\Cy2\times\Cy2}

\let\lar\longrightarrow
\let\iso\cong
\let\tensor\otimes
\let\normal\triangleleft
\let\notdiv\nmid

\def\epsn#1#2{w_{#1/#2}}
\def\eps#1#2#3{w_{#2/#3}}

\def\Kf{{M}}
\def\Kt{{M'}}

\long\def\comment#1\endcomment{}

\begin{document}

\title{Root numbers of elliptic curves in residue characteristic 2}
\author{Tim$^\dagger$ and Vladimir Dokchitser}
\date{November 29, 2006}
\subjclass[2000]{Primary 11G07; Secondary 11F80}
\thanks{$^\dagger$Royal Society University Research Fellow}
\address{Robinson College, Cambridge CB3 9AN, United Kingdom}
\email{t.dokchitser@dpmms.cam.ac.uk}
\address{Gonville \& Caius College, Cambridge CB2 1TA, United Kingdom}
\email{v.dokchitser@dpmms.cam.ac.uk}

\begin{abstract}
To determine the global root number of an elliptic curve defined over
a number field, one needs to understand all the local root numbers.
These have been classified except at places above 2, and in this paper
we attempt to complete the classification.
At places above~2, we express the local root numbers in terms of norm residue
symbols (resp. root numbers of explicit 1-dimensional characters)
in case when wild inertia acts through a cyclic (resp. quaternionic) quotient.
\end{abstract}

\phantom{}\vskip -1.4cm
\maketitle

\section{Introduction}

Let $E$ be an elliptic curve over a number field $F$.
The Hasse-Weil $L$-function $L(E/F,s)$ conjecturally satisfies a
functional equation under $s \leftrightarrow 2-s$ with the sign
given by the global root number $w(E/F)=\pm1$.
Consequently, $w(E/F)$ determines
the parity of $\ord_{s=1}L(E/F,s)$ and, assuming the Birch--Swinnerton-Dyer
conjecture, the parity of the Mordell-Weil rank of $E$ over~$F$.
This, together with the fact
that $w(E/F)$ is defined independently of any conjectures,
makes it an important invariant. 

The global root number is a product of local root numbers,
$$
  w(E/F) = \prod\nolimits_v w(E/F_v)
$$
with $w(E/F_v)$ determined by the action of the local Galois group on the
$l$-adic Tate module $T_l(E)$ for any $l\notdiv v$. To use this formula,
one needs all the local root numbers.
Their classification has been completed by Rohrlich \cite{RohG}
apart from places $v|6$ where $E$ has additive potentially good reduction,
and by Kobayashi \cite{Kob} in residue characteristic~3.
For $F_v=\Q_2, \Q_3$ all possibilities have been tabulated
by Halberstadt \cite{Hal}.

In this paper we consider the remaining case $v|2$ where $E$ has
potentially good reduction. 
It is not hard to extend
the methods of \cite{RohG,Kob} when wild inertia
acts through a cyclic quotient on $T_l(E)$
(this is always true in odd residue characteristic).
One then gets a formula for the local root number $w(E/F_v)$ in terms
of norm residue symbols (Proposition~\ref{abinertia}).

The remaining ``hard'' case is when wild inertia
acts through the quaternion group $\Qu$. In this case we express
$w(E/F_v)$ in terms of local root numbers of explicit
primitive 1-dimensional characters of cyclic extensions of degree
2,4 and 8 (Theorem \ref{main}).
Although not as satisfactory as in the other cases, this does allow one
to determine $w(E/F_v)$ and hence the global root number $w(E/F)$
for any given elliptic curve over a number field.
(It has now been implemented by the first author (T.) in Magma.)

The second author (V.) would like to thank the
Max Planck Institute for Mathematics in Bonn where some of this research
has been carried out.
We would also like to thank Mark Watkins and Steve Donnelly for pointing out
that D. Whitehouse \cite{Whi} has obtained similar results independently.

\section{Notation}

Throughout the paper $K$ is a finite extension of $\Q_2$, and $E/K$ is
an elliptic curve with additive, potentially good reduction.
For any $l\ne 2$ inertia acts on the Tate module $T_l(E)$
through a finite quotient. By \cite{ST} Cor. 2b,
$E$ has good reduction over $K(E[l])$, equivalently
the action of inertia factors through $\Gal(K(E[l])/K)$.

We will take $l=3$.
Let $L=K(E[3])$ be the field obtained by adjoining the coordinates of
the 3-torsion points of $E$, set $G=\Gal(L/K)$ and write $I \triangleleft G$
for its inertia subgroup. The action of $G$ on $E[3]$ gives an
inclusion $G\subset\GL(2,\F_3)$, under which $I\subset\SL(2,\F_3)$
from the properties of the Weil pairing.

Denote the isomorphism classes of subgroups of $\SL(2,\F_3)$
by $\Cy1$, $\Cy2$, $\Cy3$, $\Cy4$, $\Cy6$ (cyclic),
$\Qu$ (quaternion) and $\Sl$, and the remaining ones in
$\GL(2,\F_3)$ by $\Vi, \Cy8, \Di6, \Di8$ (dihedral), $\Bo (\iso\Di{12})$
(Borel), $\Sy$ (2-Sylow) and $\Gl$.

We fix an embedding $\bar\Q_3\to\C$.
Write $\zeta_8=\exp(2\pi i/8)\in\C$, $\sqrt{2}$ for the positive square
root of 2, and $\sqrt{-2}=i\sqrt{2}=\zeta_8+\zeta_8^3$.

Here is a summary of our main notation:

\smallskip
\begin{tabular}{ll}
 $K$    & finite extension of $\Q_2$\cr
 $E/K$  & elliptic curve with additive potentially good reduction\cr
 $\V$   & $=T_3(E/K)\tensor_{\Z_3}{\C}$\cr
 $L$    & $=K(E[3])$\cr
 $G$    & $=\Gal(L/K) \subset \GL(2,\F_3)$\cr
 $I$    & the inertia subgroup of $G$\cr
 $n(E)$ & valuation of the conductor of $E/K$\cr
 $\Delta, c_4, c_6$ & standard invariants of some
         fixed Weierstrass model of $E/K$\cr
 $\psi$ & $=\exp(2\pi i\tr_{K/\Q_2}(\cdot))$, an additive character $K\to\C^*$  \cr
 $n(\psi)$ & $=v(\D(K/\Q_2))$, the largest $j$ for which
             $\tr_{K/\Q_2}(\pi_K^{-j}\O_K)\subset\Z_2$\cr
 $\epsilon(\cdot)$ & local $\epsilon$-factor with respect to $\psi$ and
     a fixed Haar measure $d\mu$\cr
 $w(\cdot)$ & $= \epsilon(\cdot)/|\epsilon(\cdot)|$, the local root number
    (independent of $d\mu$)\cr
 $f_{K/\Q_2}$ & residue degree\cr
 $\mu_3\subset\bar K$ & the set of 3rd roots of unity\cr
 $\Frob_K$ & arithmetic Frobenius of $K$\cr
 $\triv{}$ & trivial representation\cr
 $(x,K'/K)\!\!\!$ & $= \pm1$, the norm residue symbol; for $K'/K$ abelian and\cr
           &$x\in K^*$ with $x^2\in N_{K'/K}{K'}$ this is 1 iff $x\in N_{K'/K}K'$ \cr
\end{tabular}

\smallskip
Finally, let $\eta: \Gal(\bar K/K)\to\C^*$ be the character that maps
the inertia subgroup to 1 (so $\eta$ is unramified) and the
arithmetic Frobenius $\Frob_K$ to $(\sqrt{-2})^{f_{K/\Q_2}}$.
The reason for this definition is the following lemma:

\begin{lemma}
\label{eta}
Suppose $G$ is non-abelian. Then the action of $\Gal(\bar K/K)$ on
$V_\eta=\V\tensor\eta^{-1}$ factors through $G$,
and $G$ acts faithfully on $V_\eta$.
Moreover,
$$
  w(E/K)=(-i)^{f_{K/\Q_2}(n(E)+2n(\psi))} w(V_\eta)\>.
$$
\end{lemma}

\begin{proof}
To begin with, $\V$ is an irreducible $\Gal(\bar K/K)$-module
($G$ is non-abelian), and the action factors through $\Gal(L^{un}/K)$.
As $\Frob_L$ is central in this group, it acts on $\V$ as
a scalar matrix $\lambda\,\text{Id}$. From
the properties of the Weil pairing, $\det(\Frob_L)=2^{f_{L/\Q_2}}$,
so $\lambda=\pm 2^{f_{L/\Q_2}/2}$.
Since $\Frob_L$ acts trivially on $E[3]$, $\lambda\equiv 1\mod 3$ so
$\lambda=(-2)^{f_{L/\Q_2}/2}$.
As $\eta$ is unramified,
$$
  \eta(\Frob_L)=\eta(\Frob_K)^{f_{L/K}}=(\sqrt{-2})^{f_{L/\Q_2}}=\lambda,
$$
so $\Gal(L^{un}/L)$ acts trivially on $\V\tensor\eta^{-1}$, and the
representation factors through $\Gal(L/K)=G$. As $G$ acts faithfully on
$E[3]$, it follows that $I$ acts faithfully on $\V$, and hence on $V_\eta$.
So the kernel of the action of $G$ on $V_\eta$ is a normal subgroup of $G$
which meets $I$ trivially. From the fact that $G/I$ is cyclic and the
classification of non-abelian subgroups of $\GL(2,\F_3)$, one sees that
this kernel is 1.

Finally, by the formula for the $\epsilon$-factor of a
twist by an unramified representation (\cite{TatN} (3.4.6)),
$$
  \epsilon(E/K)=\epsilon(V_\eta\tensor\eta)
    =\epsilon(V_\eta)^{\dim\eta}\cdot\det(\Frob_K^{-1}|\eta)^{n(E)+2n(\psi)}
$$
Now $\dim\eta=1$ and $
  \eta(\Frob_K^{-1})= (\tfrac{-i}{\sqrt 2})^{f_{K/\Q_2}}.
$
\end{proof}

\section{The structure of $\Gal(K(E[3])/K)$}

We start by classifying the possibilities
for $G=\Gal(K(E[3])/K)$ in terms of standard
invariants of $E/K$.

\begin{proposition}
\label{classg}
Depending on whether $\mu_3\subset K$, whether
$\Delta^{1/3}\in K$, and on the degrees of the irreducible factors
of $\gamma(x)=x^8-6c_4x^4-8c_6x^2-3c_4^2=\prod_i\gamma_i$ over $K$, the group
$G=Gal(L/K)$ is given by
$$
\begin{array}{|ll|ll|ll|ll|}
\hline
\multicolumn{4}{|c|}{\mu_3\subset K} & \multicolumn{4}{|c|}{\mu_3\notsubset K}\cr
\multicolumn{2}{|c}{\Delta\in K^{*3}} &
\multicolumn{2}{c|}{\Delta\notin K^{*3}} &
\multicolumn{2}{|c}{\Delta\in K^{*3}} &
\multicolumn{2}{c|}{\Delta\notin K^{*3}} \cr
\hline
(\deg\gamma_i)_i & G & (\deg\gamma_i)_i & G &(\deg\gamma_i)_i & G &(\deg\gamma_i)_i & G \cr
\hline
    (2,2,2,2)   & \Cy2 & (1,1,3,3)   & \Cy3 &
    (2,2,4)     & \Cy2\times\Cy2 & (1,1,6)     & \Di6           \cr
    (4,4)       & \Cy4 & (2,6)       & \Cy6 &
    (4,4)       & \Di8           & (2,3,3)     & \Di6           \cr
    (8)         & \Qu  & (8)         & \Sl  &
    (8)         & \Cy8\>\,\text{or}\>\,\Sy & (2,6)       & \Bo            \cr
    &&&&&& (8)         & \Gl            \cr
\hline
\end{array}
$$
\end{proposition}

\begin{proof}

From the non-degeneracy of the Weil pairing on $E[3]$ it follows that
$G\subset\SL(2,\F_3) $ iff $\mu_3\subset K$. Also
$\Delta^{1/3}\in L$, and $\Delta^{1/3}\in K$ iff $3\notdiv |G|$
(look at the resolvent cubic of the 3-torsion polynomial, or see \cite{Kra}).

A nonzero 3-torsion point $P=(x_3,y_3)\in E(\bar K)[3]$ is an
inflection point on $E$ and the slope $\kappa$ of the tangent line to $P$
lies in $K(x_3,y_3)$. Conversely, it is easy to see that $x_3,y_3\in K(\kappa)$,
and that $2\kappa$ is a root of $\gamma(x)$ (see e.g. \cite{TD},
proof of Thm. 1 with $a=-c_4/48$ and $b=-c_6/864$).
It follows that $L$ is the splitting field of $\gamma(x)$ over $K$ and
the action of $G$ on its roots is the same as that on $E[3]\setminus\{O\}$.

Now the statement follows by inspection of the possible subgroups
of $\GL(2,\F_3)$ and lengths of their orbits under their
natural action on $\F_3^2$.
\end{proof}

There is one case ($G=\Cy8$ or $G=\Sy$), when the proposition
does not determine $G$ uniquely. To deal with it, we have

\begin{lemma}
\label{classg2}
Assume that $x^3-12^3\Delta$ has exactly one root $\delta$ in $K$
and that $\gamma(x)$ is irreducible. Then the following conditions are equivalent
\begin{enumerate}
\item $G=\Cy8$,
\item $\gamma(x)$ is reducible over $K(\mu_3)$,
\item The numbers $-3(c_4-\delta)$ and $-3(c_4^2+c_4\delta+\delta^2)$ are squares
in $K$ (if $c_6\ne0$, the two conditions are equivalent).
\end{enumerate}
\end{lemma}

\begin{proof}
We have $\mu_3\notsubset K$ and $\Delta^{1/3}\in K$, so $G\in\{\Cy8,\Sy\}$ by
Proposition~\ref{classg}.
The group $\Gal(L/K(\mu_3))$ is $\Cy4$ for $G=\Cy8$ and $\Qu$ for $G=\Sy$, so
the equivalence $(1)\iff (2)$ follows from Proposition \ref{classg} over $K(\mu_3)$.

The equivalence $(1)\iff (3)$ follows from the classification of Kraus (\cite{Kra}, Thm. 3).
\end{proof}

\section{Abelian inertia}

The case of abelian inertia is essentially due to Rohrlich, and the arguments
below are also used by Kobayashi in residue characteristic 3. (Note the misprints $2\xi$ instead of
$\xi$ in Prop.~5.2, and $a/2+v(2)$ instead of $a/2$ in Prop.~6.1 of
\cite{Kob}.)

\begin{proposition}
\label{abinertia}
Assume that $G$ is either abelian, dihedral, or 
$L/K$ is not totally ramified and $G=\Qu$. Then
\begin{itemize}
\item[(a)] The inertia subgroup $I$ is cyclic.
\item[(b)] If $G$ is abelian, then $w(E/K)=(-1,L/K)$.
\item[(c)] If $G$ is non-abelian, write the unique quadratic unramified
extension of $K$ as $K(\xi)$ with $\xi^2\in\O_K^*$. Then
$w(E/K)=(-1)^{n(E)/2}(\xi,L/K(\xi))$.
\end{itemize}

\begin{proof}
(a) This follows from the fact that $I\subset\Sl$ and the classification of
subgroups of $\Sl$.

(b) Since $G$ is abelian, the action of $\Gal(\bar K/K)$ on $\V$
factors through an abelian quotient. As $\V$ is semisimple,
it must decompose as $\chi\oplus\chi^{-1}||\cdot||_K$,
with $||\cdot||_K$ the cyclotomic character.
In this case, determinant formula (\cite{RohE}, p.145) implies
$w(E/K)=\chi(-1)=(-1,L/K)$.

(c) Here $G$ is either dihedral
($\Di6$, $\Di8$ or $\Bo$) or $\Qu$ and $I$
must be cyclic of index 2 in $G$, since
these groups have no cyclic quotients of higher order and
$I\subset\SL(2,\F_3)$ so $I\ne G$.

By Lemma \ref{eta}, $V_\eta=\V\tensor\eta^{-1}$
is the 2-dimensional faithful irreducible representation of $G$.
Denote by $\nu$ the quadratic unramified character of $G$,
and by $\phi$ a primitive character of $I$. Then $V_\eta=\Ind\phi$ and,
by inductivity in degree 0,
$$
  w(\phi)
    = \frac{w(\phi)}{w(\triv{K(\xi)})}
    = \frac{w(V_\eta)}{w(\triv K)w(\nu)}
    = \frac{w(V_\eta)}{\nu(\Frob_K^{-1})^{n(\psi)}} = (-1)^{n(\psi)} w(V_\eta)\>.
$$
So, by Lemma \ref{eta},
$$
\begin{array}{rcl}
  w(E/K)&=&(-i)^{f_{K/\Q_2}(n(E)+2n(\psi))} (-1)^{n(\psi)} w(\phi)\cr
        &=&(-1)^{(f_{K/\Q_2}+1)n(\psi)} (-1)^{f_{K/\Q_2}n(E)/2} w(\phi).
\end{array}
$$

Suppose $G$ is dihedral. Since $G\not\subset\SL(2,\F_3)$, the residue
degree $f_{K/\Q_2}$ is odd.
By a theorem of Fr\"ohlich-Queyrut (\cite{FQ}, Lemma 1 and Thm. 3),
$$
  w(\phi)=\phi(\theta_{K(\xi)}(\xi))=(\xi,L/K(\xi))\>,
$$
where $\theta_{K(\xi)}$ is the local reciprocity map on $K(\xi)^*$.
Combining these formulae, $w(E/K)=(-1)^{n(E)/2}(\xi,L/K(\xi))$.

Finally, if $G=\Qu$, then $\mu_3\subset K$, and
$w(V_\eta)=w(\V)$.
Moreover, the theorem of Fr\"ohlich-Queyrut does not apply directly,
since $(\phi\circ\theta_{K(\xi)})|_{K^*}=\nu\circ\theta_K$, which is
non-trivial.
In this case, we twist $\phi$ by an unramified quadratic character of $K(\xi)^*$,
and a computation as in \cite{RohV} proof of Prop. 2 (esp. p.130)
yields the result in this case.
\end{proof}

\end{proposition}

\begin{remark}
In some cases, this result can be further simplified.
For $G=C_3$, $w(E/K)=1$, for $G=D_6$, $w(E/K)=-1$, and for
$G=\Bo$, $w(E/K)=(-1,K(\kappa)/K)$ with $\kappa$ a root of a the quadratic
factor of $\gamma(x)$ as in Proposition \ref{classg}
(see \cite{TV-P}, proof of Lemma 10).
\end{remark}

\section{Non-abelian inertia}

From now on, suppose that we are in one of the cases not covered
by Proposition \ref{abinertia}, so
$G=\Sy,\Sl,\Gl$ or $G=I=\Qu$. Since $I\subset\SL(2,\F_3)$ and
$G/I$ is cyclic, it follows that $\Qu\subset I$, in particular $I$
is non-abelian.

\begin{lemma}
\label{quext}
Assuming $\Qu\subset I$,
the extension $L/K$ is built up from the
following subfields (we choose the roots once and for all):
$$
\begin{array}{ll}
  \>L\!=\!K_4(\mu_3)\cr
  \>\>|\cr
  K_4\!=\!K_3(y_3),& y_3=\sqrt{x_3^3/3-c_4x-2c_6/3},\>x_3=\sqrt{A}+\sqrt{C}\cr
  \llap{$\scriptstyle 2\!\!$}\>\>|\cr
  K_3\!=\!K_2(\sqrt{C}),& C=2c_4+12 \Delta^{1/3}+2\sqrt{B}\cr
  \llap{$\scriptstyle 2\!\!$}\>\>|\cr
  K_2\!=\!K_1(\sqrt{A},\sqrt{B}), &
    A=c_4\!-\!12 \Delta^{1/3}, B=c_4^2\!+\!12 c_4 \Delta^{1/3}\!+\!(12\Delta^{1/3})^2 \cr
  \llap{$\scriptstyle 2\!\!$}\>\>|&AB=c_6^2\cr
  K_1\!=\!K(\Delta^{1/3})\cr
  \>\>|\cr
  \,K
\end{array}
$$
The extensions $K_4/K_3$, $K_3/K_2$ and $K_2/K_1$
are ramified quadratic,
$K_1/K$ is either trivial $(3\notdiv \#G)$, unramified $(3|v(\Delta))$ or tamely
ramified $(3\notdiv v(\Delta))$, and $L/K_4$ is either trivial $(\mu_3\subset K)$ or
unramified quadratic $(\mu_3\notsubset K)$.
\end{lemma}

\begin{proof}
Essentially everything here is due to Kraus \cite{Kra}, Thm. 3. He defines
$A, B$ and $C$, and exhibits their relation to the coordinates of the
3-torsion points (\cite{Kra}, pp. 371, 373). Since $L/K_4$ and $K_1/K$
cannot contribute to the 2-part of $I$, it follows that
$K_4/K_3$, $K_3/K_2$ and $K_2/K_1$ are ramified, in consistence with
\cite{Kra}, Thm. 3.
\end{proof}

\begin{theorem}
\label{main}
Assume that $\Qu\subset I$, and keep the notation of Lemma \ref{quext}.
Set $\Kf=K_1(\sqrt{-3A},\sqrt{-3B})$. Then
$$
  w(E/K)=(-i)^{f_{K/\Q_2}(n(E)+2n(\psi))}
  \>\cdot\> \epsn{L}{\Kf} \>\cdot\>
  \left\{
    \begin{array}{ll}
      \epsn{K(x_3,y_3)}{K(x_3)}^{-1}, & I=\Sl \cr
      \epsn{\Kf}{K_1}, & G\ne\Gl \cr
    \end{array}
  \right.
$$
Here $\epsn{K(x_3,y_3)}{K(x_3)}$ and $\epsn{\Kf}{K_1}$ denote the
root numbers of the non-trivial character of the corresponding
ramified quadratic extensions, and
$\epsn{L}{\Kf}$
is the root number of a character $\chi$ of the cyclic
extension $L/\Kf$ with the property that
$\Ind_{\Gal(L/\Kf)}^{\Gal(L/K)}\chi$ contains a copy of $V_\eta$.
If $\mu_3\subset K$, then any primitive $\chi$ has this property.
\end{theorem}

\begin{proof}
By Lemma \ref{eta}, $w(E/K)=(-i)^{f_{K/\Q_2}(n(E)+2n(\psi))}w(V_\eta)$,
so it remains
to show that $w(V_\eta)$ is given by one of the two products
of root numbers in the statement of the theorem.
Note here that when $I=G=\Sl$, both formulae are valid.

The first formula is covered by Lemma \ref{lemgl} and Lemma \ref{lemsl}.

Now we prove the second formula for $G\ne\Gl$.
If $I=\Qu$ and $G\in\{\Qu,\Sy\}$, this is Lemma \ref{lemsy}.
In the remaining cases $G=\Sl$, so $\mu_3\subset K$ and
$K_1=K(\Delta^{1/3})$ is abelian over $K$ of degree 3.
We know that the formula holds for $E/K_1$ by the $I=G=\Qu$ case.
On the other hand,
$2|f_{K/\Q}$, $n(E/K)\equiv n(E/K_1)\mod 2$, and
$w(E/K)=w(E/K_1)$ by a result of Kramer--Tunnell
\cite{KT}, proof of Prop. 3.4.
\end{proof}


\begin{lemma}
\label{lemgl}
If $I=\Sl$ and $\mu_3\notsubset K$ (equivalently $G=\Gl$), then
$$
  w(V_\eta)=\epsn{L}{\Kf}/\epsn{K(x_3,y_3)}{K(x_3)}.
$$
\end{lemma}

\begin{proof}
From the character table of $G=\GL(2,\F_3)$ (Table \ref{glchartable}),
$V_\eta$ is either $\rho'_4$ or $\rho''_4$,
these two being the only two irreducible faithful representations of $G$
of dimension 2.
To write $V_\eta$ as a combination of inductions of 1-dimensional
characters, it is necessary to use a subgroup $C_8\subset G$, because
induced characters from all other subgroups do not distinguish between
$\rho'_4$ and $\rho''_4$. Consider the 2-Sylow subgroup of $G$
$$
  \Sy=\Gal(L/K_1)\subset\Gal(L/K)=G.
$$
It has 3 index 2 subgroups ($\Cy8, \Di8$ and $\Qu$), that correspond to 3
quadratic extensions of $K_1$, namely $K_1(\mu_3), K_2$ and $\Kf$.
Clearly $\Qu\subset\Sl\subset G$ corresponds to $K_1(\mu_3)$, and
Kraus' formulae (\cite{Kra} p. 373) show that
$\Gal(L/K_2)$ does not act transitively on the roots of the 3-torsion
polynomial, hence it is dihedral by Proposition \ref{classg}.
Thus $\Gal(L/\Kf)\iso\Cy8$.

\begin{table}[h]
$$
\begin{array}{c|rrrrrrrrrr}
\text{order}  & 1 &  2 &  2 &  3 &  4 &  6 &  8 & 8    \cr
\text{\#elts} & 1 &  1 & 12 &  8 &  6 &  8 &  6 & 6    \cr
\hline
\rho_1        & 1 &  1 &  1 &  1 &  1 &  1 &  1 & 1    \cr
\rho_2        & 1 &  1 & -1 &  1 &  1 &  1 & -1 & -1   \cr
\rho_3        & 2 &  2 &  0 & -1 &  2 & -1 &  0 & 0    \cr
\rho'_4       & 2 & -2 &  0 & -1 &  0 &  1 & \sqrt{-2} & -\sqrt{-2} \cr
\rho''_4      & 2 & -2 &  0 & -1 &  0 &  1 & -\sqrt{-2} & \sqrt{-2} \cr
\rho_5        & 3 &  3 &  1 &  0 & -1 &  0 & -1 & -1   \cr
\rho_6        & 3 &  3 & -1 &  0 & -1 &  0 &  1 & 1    \cr
\rho_7        & 4 & -4 & 0 &   1 &  0 & -1 &  0 & 0    \cr
\end{array}
$$
\caption{Character table of $\GL(2,\F_3)$}
\label{glchartable}
\end{table}

If $\chi$ is a primitive character of $\Cy8$, then
$\Ind_{\Cy8}^G\chi\iso\rho_4\oplus\rho_7$ with $\rho_4\in\{\rho'_4,\rho''_4\}$.
Replacing $\chi$ by $\chi^{-1}$ if necessary (this changes the induction),
assume that $V_\eta$ is contained in $\Ind_{\Cy8}^G\chi$.

If $\Bo\subset G$ is the upper-triangular matrices, let $\det_{\Bo}$ be the determinant
character, and $\sigma_{\Bo}$ the ``top left corner'' character, both of
order 2. Let $\tau_{\Sl}$ be a non-trivial 1-dimensional character
of $\Sl\subset G$. A computation with inductions shows that
$$
  \Ind_{\Cy8}^{G} ( \chi\ominus \triv{\Cy8} ) \oplus
  \Ind_{\Bo}^{G} ( \det\nolimits_{\Bo}\ominus\sigma_{\Bo} ) \oplus
  \Ind_{\Sl}^{G} ( \tau_{\Sl}\ominus \triv{\Sl} ) \iso
$$
$$
  \iso (V_\eta\oplus\rho_7\ominus\rho_1\ominus\rho_3\ominus\rho_6) \oplus
  (\rho_2\oplus\rho_6\ominus\rho_7) \oplus
  (\rho_3\ominus\rho_1\ominus\rho_2) = V_\eta \ominus 2\rho_1\>,
$$
hence
$$
  w(V_\eta) =
     \frac{w(\chi)}{1}\cdot
     \frac{w(\det\nolimits_{\Bo})}{w(\sigma_{\Bo})}\cdot
     \frac{w(\tau_{\Sl})}{1}.
$$
By definition, $\eps8{L}{\Kf}=w(\chi)$.
Since $\Sl=\Gal(L,K(\mu_3))$ and $K(\mu_3,\Delta^{1/3})/K$ is dihedral,
$w(\tau_{\Sl})=1$ by Fr\"ohlich-Queyrut's Theorem
(\cite{FQ}, Thm. 3 and proof of Lemma 1).

Now, $\Bo=\Gal(L,K(x_3))$ is the subgroup leaving one 3-torsion subgroup
invariant, so $w(\sigma_{\Bo})=\eps2{K(x_3,y_3)}{K(x_3)}$. Finally,
$$
\begin{array}{rcl}
  w(\det\nolimits_{\Bo}) &\!\!=\!\!&
  (-1)^{n(\psi_{K(x_3)})} = (-1)^{v(\D(K(x_3),\Q_2))} \cr
  &\!\!=\!\!&(-1)^{v(\D(K(x_3),K))+[K(x_3):K]v(\D(K,\Q_2))} = (-1)^{v(\D(K(x_3),K))} = 1.
\end{array}
$$
The last equality uses that the parity of $v(\D(K(x_3),K))$ is the same
as that of the $K$-valuation of the discriminant of the 3-torsion polynomial
$P(t)$, which defines $K(x_3)/K$; but $\Delta(P(t))=-2^{12}3^3\Delta(E)^2$
has even valuation over any 2-adic field.
\end{proof}


\begin{lemma}
\label{lemsl}
If $I=\Sl$ and $\mu_3\subset K$ (so $G=\Sl$), then
$$
  w(V_\eta)=\epsn{L}{\Kf}/\epsn{K(x_3,y_3)}{K(x_3)}.
$$
\end{lemma}

\begin{proof}
From the character table of $G=\SL(2,\F_3)$ (Table \ref{slchartable}),
we have $V_\eta\iso\rho_4$, since in the other two-dimensional irreducible
representations $\rho_5$ and $\rho_6$ of $G$ elements of order 3 do not have
determinant 1 (look at the traces).

\begin{table}[h]
$$
\begin{array}{c|rrrrrrrrrr}
\text{order}  & 1 & 2 &  3 &  3 & 4 &  6 &  6\cr
\text{\#elts} & 1 & 1 &  4 &  4 & 6 &  4 &  4\cr
\hline
\rho_1 & 1 &  1 &        1 & 1 & 1 & 1 & 1\cr
\rho_2 & 1 &  1 & \omega^2 & \omega & 1 & \omega & \omega^2\cr
\rho_3 & 1 &  1 & \omega   & \omega^2 & 1 & \omega^2 & \omega\cr
\rho_4 & 2 & -2 &       -1 & -1 & 0 & 1 & 1\cr
\rho_5 & 2 & -2 & -\omega^2 & -\omega & 0 & \omega & \omega^2\cr
\rho_6 & 2 & -2 & -\omega & -\omega^2 & 0 & \omega^2 & \omega\cr
\rho_7 & 3 &  3 & 0 & 0 & -1 & 0 & 0\cr
\end{array}
$$
\caption{Character table of $\SL(2,\F_3)$}
\label{slchartable}
\end{table}

Take $\Cy4\subset \Qu \normal G$, and let $\chi$ be one of the two of its order
4 characters. Consider also a Borel subgroup $\Cy6$ of $G$
and its character $\sigma$ of order 2. Then one computes
$$
\begin{array}{l}
  \Ind_{\Cy4}^{G} ( \chi\ominus \triv{\Cy4} ) \oplus
  \Ind_{\Cy6}^{G} ( \triv{\Cy6}\ominus \sigma ) \cr
  \iso (\rho_4\oplus\rho_5\oplus\rho_6) \ominus (\rho_1\oplus\rho_2\oplus\rho_3\oplus\rho_7)
    \oplus
    (\rho_1\oplus\rho_7) \ominus (\rho_5\oplus\rho_6) \cr
  \iso V_\eta \ominus (\rho_2\oplus\rho_3) \>.
\end{array}
$$
Thus,
$$
  w(V_\eta) =
     w(\chi) w(\sigma)^{-1} w(\rho_2\oplus\rho_2^{-1})\>.
$$
By the determinant formula,
$w(\rho_2\oplus\rho_2^{-1})=\rho_2(-1)=1$, as $\rho_2(-1)$ is a cube
root of unity that squares to 1. Finally, note that
$$
  w(\chi)=\epsn{L}{\Kf} \qquad
  \text{and}\qquad w(\sigma)=\epsn{K(x_3,y_3)}{K(x_3)} \>.
$$
\end{proof}

\begin{lemma}
\label{lemsy}
If $I=\Qu$ and $G\in\{\Qu,\Sy\}$, then
$$
  w(V_\eta)=
  \epsn{L}{\Kf}\,\epsn{\Kf}{K_1}.
$$
\end{lemma}
\begin{proof}
If $G=\Sy$ and $\chi$ is a primitive
character of $\Gal(L/\Kf)=\Cy8$, then $\Ind_{\Cy8}^G\chi$ is isomorphic
either to $\rho_4'|_{\Sy}$ or to $\rho_4''|_{\Sy}$ from Table \ref{glchartable}.
As in the $G=\Gl$ case, these are the only two faithful 2-dimensional
irreducible representations of $\Sy$,
and we can choose $\chi$ with $\Ind_{\Cy8}^G\chi=V_\eta$. Now,
$$
  \frac{w(V_\eta)}{w(\Ind_{\Cy8}^G \triv{\Cy8})} =
  \frac{w(\chi)}{w(\triv{\Cy8})},
$$
implying the asserted formula.
The case $G=\Qu$ is analogous.
\end{proof}

\section{Determining which character to take}

Assume that $G=\Sy$ or $G=\Gl$, so $[G:I]=2$. Then $\Gal(L/\Kf)\iso\Cy8$,
which has 4 characters of order 8,
say $\chi_1, \chi_1^3, \chi_1^{-3}, \chi_1^{-1}$.
We have (cf. Table~\ref{glchartable})
$$
  \Ind_{\Cy8}^G \chi_1 = \Ind_{\Cy8}^G \chi_1^3 = \rho_4'\oplus \sigma, \qquad
  \Ind_{\Cy8}^G \chi_1^{-3} = \Ind_{\Cy8}^G \chi_1^{-1} = \rho_4''\oplus \sigma,
$$
with $\sigma=0$ for $G=\Sy$ and $\sigma=\rho_7$ for $G=\Gl$. So two of the characters
have the property that their induction contains $V_\eta$.
It may in fact happen that $w(\chi^{-1})=-w(\chi)$,
so it is important that we choose the right character when computing
the root number. Pick one of the
4 characters $\chi$, and compose it with local reciprocity,
$$
  \chi_\Kf = \chi\circ\theta_\Kf: \Kf^* \lar \C^* \>.
$$
Then we have the following criterion:

\begin{proposition}
Let $\Kt/\Kf$ be a quartic totally ramified extension with $\Kt\subset L^{un}$,
let $\pi_{\Kt}$ be its uniformiser, let $\pi_\Kf=N_{\Kt/\Kf}\pi_{\Kt}$ and
define $\delta=\pm 1$ by
$$
  \chi_\Kf(\pi_\Kf^{-1})+\chi_\Kf(\pi_\Kf^{-1})^3=\delta \sqrt{-2}.
$$
$E$ has good reduction over $\Kt$, and let $\tilde E/\F_q$ denote the reduced
curve. The following conditions are equivalent:
\begin{enumerate}
\item $\Ind_{\Cy8}^G\chi$ contains $V_\eta$.
\item The arithmetic Frobenius $\theta_\Kf(\pi_\Kf^{-1})\in\Gal(L^{un}/\Kt)$
acts with trace $\delta$ (mod 3) on $E[3]$.
\item $|\tilde E(\F_q)|=q+1-\delta(-2)^{(1+f_{K/\Q_2})/2}$.
\end{enumerate}
\end{proposition}

\begin{proof}
Suppose $\Ind_{\Cy8}^G\chi$ contains $V_\eta$, so
$\Res_{\Kf/K}V_\eta=\chi\oplus\chi^3$. Then with $\chi'=\Res_{\Kt/\Kf}\chi$
and $\chi_{\Kt}$ for the corresponding character on ${\Kt}^*$,
$$
\begin{array}{rcl}
  \delta\sqrt{-2} &=&
  \chi_\Kf(\pi_\Kf^{-1})+\chi_\Kf(\pi_\Kf^{-1})^3\cr
   &=& \chi_{\Kt}(\pi_{\Kt}^{-1})+\chi_{\Kt}(\pi_{\Kt}^{-1})^3 \cr
   &=&\chi'(\Frob_{\Kt}) + (\chi'(\Frob_{\Kt}))^3\cr
   &=& \tr\bigl(\Frob_{\Kt}\!\bigm|\!\Res_{\Kt/K} (V_\eta)\bigr)\cr
   &=& (\Res_{\Kt/K}\eta^{-1})(\Frob_{\Kt}) \cdot
       \tr\bigl(\Frob_{\Kt}\!\bigm|\!\Res_{\Kt/K} V\bigr) \cr
   &=& a (\sqrt{-2})^{-f_{\Kt/\Q_2}}\cr
   &=& a (\sqrt{-2})^{-f_{K/\Q_2}}.
\end{array}
$$
In other words, if $\Ind_{\Cy8}^G\chi$ contains $V_\eta$, then
$a=\delta(-2)^{(1+f_{K/\Q_2})/2}$.
Otherwise this formula holds for $\chi^{-1}$ (for which $\delta$ changes
sign), so $(1)\iff(2)\iff(3)$.
\end{proof}


\begin{thebibliography}{10}


\bibitem{TV-P}
T. Dokchitser, V. Dokchitser,
Parity of ranks for elliptic curves with a cyclic isogeny, preprint May 2006,
arxiv: math.NT/0604149.

\bibitem{TD}
T. Dokchitser, Ranks of elliptic curves in cubic extensions, to appear in
Acta Arith.

\bibitem{FQ}
A. Fr\"ohlich, J. Queyrut, On the functional equation of the Artin
$L$-function for characters of real representations,
Invent. Math. 20 (1973), 125--138.

\bibitem{Hal}
E. Halberstadt, Signes locaux des courbes elliptiques en 2 et 3,
C. R. Acad. Sci. Paris S\'erie I Math. 326 (1998), no. 9, 1047--1052.

\bibitem{Kob}
S. Kobayashi, The local root number of elliptic curves with wild
ramification, Math. Ann. 323 (2002), 609--623.

\bibitem{KT}
K. Kramer, J. Tunnell, Elliptic curves and local $\epsilon$-factors,
Compositio Math. 46 (1982), 307--352.

\bibitem{Kra}
A. Kraus, Sur le d\'efaut de semi-stabilit\'e des courbes elliptiques
\`a r\'eduction additive, Manuscripta Math. 69 (1990), no. 4, 353--385.

\bibitem{RohE}
D. Rohrlich, Elliptic curves and the Weil-Deligne group, in Elliptic curves
and related topics, 125--157, CRM Proc. Lecture Notes 4, Amer. Math. Soc.,
Providence, RI, 1994.

\bibitem{RohG}
D. Rohrlich, Galois Theory, elliptic curves, and root numbers, Compositio
Math. 100 (1996), 311--349.

\bibitem{RohV}
D. Rohrlich, Variation of the root number in families of elliptic
curves, Compositio Math. 87 (1993), 119--151.

\bibitem{ST}
J.-P. Serre, J. Tate, Good reduction of abelian varieties,
Annals of Math. 88 (1968), 492--517.

\bibitem{TatN}
J. Tate, Number theoretic background, Proceedings of Symposia
in Pure Mathematics 33 (1979), Part 2, 3--26.

\bibitem{Whi}
D. Whitehouse, Root numbers of elliptic curves over 2-adic fields, preprint.

\end{thebibliography}
\end{document}